\documentclass[a4paper,fleqn]{cas-sc}

\usepackage[authoryear]{natbib}

\usepackage{graphicx}
\usepackage{subcaption}
\usepackage{multirow}
\usepackage{tikz}
\usetikzlibrary{arrows.meta, positioning}

\definecolor{ForestGreen}{RGB}{34,139,34}

\usepackage{etoolbox}
\AtBeginDocument{%
  \setlength{\abovedisplayskip}{4pt}
  \setlength{\belowdisplayskip}{4pt}
  \setlength{\abovedisplayshortskip}{4pt}
  \setlength{\belowdisplayshortskip}{4pt}
}

\makeatletter
\def\fps@figure{!htbp}
\makeatother
\usepackage{placeins}

\begin{document}
\let\WriteBookmarks\relax
\let\printorcid\relax   

\shorttitle{User-driven framework for EV charging}

\shortauthors{Zhou et~al.}

\title [mode = title]{A user-driven pricing and scheduling framework for public electric vehicle charging}

\author[1]{Fangting Zhou}

\author[1]{Jiaming Wu}

\affiliation[1]{organization={Architecture and Civil Engineering, Chalmers University of Technology},
    city={Gothenburg},
    country={Sweden}}

\author[2]{Balázs Kulcsár}

\affiliation[2]{organization={Electrical Engineering, Chalmers University of Technology},
    city={Gothenburg},
    country={Sweden}}

\begin{abstract}
Public electric vehicle (EV) charging infrastructure has expanded rapidly, {yet utilization across charging stations remains uneven and often inefficient.}
Existing operator-determined pricing schemes offer limited flexibility to coordinate heterogeneous user demand within constrained capacity. This study proposes a user-driven pricing and scheduling framework {for public EV charging.}
{Users submit advance bids specifying acceptable time slots, bid prices, and quantity bounds. Based on these bids, the charging operator determines prices and charging slot assignments. After observing the outcomes, users decide whether to accept the resulting prices and allocations.
The operator’s decision problem incorporates profit objectives, user participation requirements, and capacity constraints across charging levels and time slots.}
The framework captures a three-stage interaction {involving user bidding, operator decisions, and user acceptance.}
Numerical case studies reveal trade-offs among user acceptance, operator revenue, and charging capacity utilization under different charging levels.
{The findings provide practical guidance and insights into designing flexible pricing schemes that better accommodate heterogeneous user preferences while improving system efficiency.}
\end{abstract}


\maketitle

\section{Introduction}

Electric vehicle adoption is accelerating rapidly worldwide. In 2024, over 17 million electric vehicles were sold worldwide, pushing the electrified share of new car sales above 20\% \citep{IEA2025}. To support this rapid expansion, charging facilities has been scaling up to relieve range anxiety. Although home and depot charging remain dominant due to their affordability, public charging facilities have also grown explosively to further enable the transition. Between 2022 and 2024, the global stock of public chargers doubled to more than 5 million units \citep{IEA2025}. Even among EV owners with home access, the survey suggests that roughly half of charging occurs away from home, highlighting the strategic importance of public charging networks \citep{ChargeLab2024Survey}.

A common motivation for purchasing an electric vehicle is its lower operating cost compared to conventional cars \citep{Hardman2018}, making charging cost a central factor shaping user behavior and infrastructure usage.
Specifically, pricing has a significant influence on when, where, and how much drivers charge, thereby affecting the spatial and temporal utilization of charging infrastructure \citep{visaria2022, zhou2025}. Pricing also feeds back into power system operations by shifting demand in time, managing load peaks, and influencing capacity requirements \citep{visaria2022}. As a result, pricing is not merely an operational consideration, but a strategic factor for coordinating charging demand and supporting transport–energy system integration \citep{shao2025}.

Despite its strategic importance, public charging in practice often faces challenges related to cost, predictability, and system-wide efficiency \citep{letmathe2025,zhou2025}. In many regions, public charging remains less attractive or less predictable than home charging, resulting in low utilization of public chargers \citep{baresch2019,cui2024}.
At the same time, empirical studies suggest that public charging is still operated in a largely passive and uncoordinated manner: station utilization is {largely} determined by users’ arrival patterns and dwell times, with limited mechanisms to actively coordinate charging demand \citep{zhou2025}. This reactive operation model contributes to congestion during peak periods, underutilization during off-peak times, and overall inefficiencies in infrastructure use. {From a system perspective, the key challenge lies not only in expanding charging infrastructure but also in better coordinating when and how users charge within existing charging capacity.}
Simply adding more chargers {is unlikely to} fully resolve these issues. 
This highlights the need for more flexible pricing and coordination mechanisms that can better match user preferences with available charging resources.
{However, effectively integrating user preferences with system-level constraints in charging operations remains an open challenge.}

Recent advances in smart infrastructure and digital mobility platforms have enabled more interactive EV charging strategies \citep{shahriar2020}. While traditional charging systems rely on operator-defined pricing schemes, recent studies have explored user-side participation through game-theoretic bidding or auction-based coordination mechanisms \citep{cao2018,zhang2020,gao2021,kim2022,kim2023}. These methods allow users to express flexibility or respond to incentives but often assume single-stage interaction or overlook system-level coordination.

To address these gaps, we propose a bilateral negotiation framework that integrates personalized user bids with centralized operator decisions under system-level constraints. Unlike auction-based mechanisms or fixed-price scheduling, this approach enables structured interaction: each user submits a bid specifying a proposed price, desired energy quantity, and acceptable time window.
{The operator then determines whether to accept, reject, or counter the bid while accounting for time- and rate-dependent operating costs and capacity constraints.}
The user subsequently evaluates the counteroffer and decides based on individual utility.

This new paradigm captures the decentralized nature of user behavior while enabling explicit coordination between supply and demand under system-level constraints.
{It explicitly incorporates user flexibility into charging decisions, }supports personalized demand shaping, and enables more coordinated and balanced utilization of charging resources. From a modeling perspective, it introduces a new class of bilateral decision processes that integrate user bids with centralized scheduling in EV charging infrastructure planning.

The \textbf{contributions} of this paper are threefold:
\begin{itemize}
    \setlength{\itemsep}{0pt}
    \setlength{\parsep}{0pt}
    \setlength{\parskip}{0pt}
    \item 
    {Propose a user-driven framework that enables personalized bidding and coordinated pricing and scheduling decisions for public EV charging;}
    \item 
    {Develop a system-level decision model that integrates user-submitted bids with centralized operator decisions under infrastructure and operational constraints;}
    \item 
    {Provide analytical insights into how user participation and coordinated pricing can improve charging infrastructure utilization and allocation efficiency.}
\end{itemize}

The remainder of the paper is organized as follows. Section \ref{Section2} reviews related work on the scheduling and pricing of EV charging, with a focus on user-side incentives and two-sided coordination mechanisms. Section \ref{Section3} presents the methodology, including the user bidding process, the operator’s optimization problem, and the user’s response and acceptance behavior. Next, Section \ref{Section4} describes the experimental setup and numerical scenarios, covering both homogeneous and heterogeneous user settings under single and multiple charging rate configurations. A larger-scale case study is also included to illustrate system-level behaviors. Finally, Section \ref{Section5} concludes the paper and outlines the directions for future research.

\section{Literature review} \label{Section2}

Prior research on EV charging coordination can be broadly grouped into three categories: (i) centralized or operator-controlled scheduling and pricing mechanisms (hereafter referred to as operator-driven), (ii) market-based and preference-driven user models, and (iii) two-sided or user-driven mechanisms. Operator-driven approaches typically emphasize centralized scheduling and predefined pricing, while market- and preference-driven studies focus on user heterogeneity and behavioral responses. Two-sided or user-driven mechanisms investigate negotiations or structured interactions between users and operators. This categorization provides a structured view of the literature, and we conclude the section with a summary that highlights remaining research gaps.

\subsection*{Operator-driven scheduling and pricing} 

Most existing studies adopt a centralized or aggregator-controlled perspective, where prices and charging schedules are predefined by operators and users respond passively. A recent review confirms that the majority of charging strategies remain centralized and operator-controlled \citep{christensen2025}.
These include time-of-use pricing, real-time pricing, and incentive-based load shifting approaches \citep{wu2020,wang2019,adetunji2023,zhou2025collaborative}. Although effective in certain contexts, these methods largely assume passive user behavior and do not account for personalized flexibility or willingness to pay.

To improve scheduling efficiency, recent studies have incorporated advanced optimization and machine learning methods into operator-driven pricing. For instance, \cite{ren2023} combined a Long Short-Term Memory network with an improved linear programming model to design dynamic electricity prices and optimize V2G scheduling, reducing user costs and balancing grid load. \cite{cedillo2022} proposed a bi-level optimization framework for solar-powered charging stations, where the upper level determines dynamic tariffs and energy bids, while the lower level computes EV schedules to enable participation in balancing services. \cite{wu2023} investigated infrastructure-constrained charging, explicitly considering the limited number of chargers and station power limits, and employed an adaptive genetic algorithm to assign charging piles and minimize costs.

In addition, reinforcement learning has been widely applied to dynamic pricing problems. \cite{zhao2021} modeled the dynamic pricing process as a Markov decision problem and applied deep reinforcement learning to adjust charging prices in real time, with quality-of-service metrics such as queuing behavior taken into account. Similarly, \cite{zulfiqar2024} introduced an enhanced multi-agent neural network trained with reinforcement learning to optimize EV scheduling under dynamic pricing, achieving both cost savings and improved operator profits.

Prior studies have also recognized that EV users are sensitive to charging prices. For example, \cite{zhang2018} analyzes the relationship between price levels, service drop rates, and the efficiency of station operations, and designs a pricing scheme to minimize user dropout. These works motivate the need to incorporate user-side price responsiveness into charging strategies.
However, even when user responsiveness is acknowledged (e.g., \cite{deng2022}), most studies still assume operator-defined pricing and lack mechanisms for direct user bidding or negotiation.

\subsection*{Market-based and preference-driven user models}

In contrast to purely operator-controlled strategies, some studies explicitly account for user heterogeneity and price sensitivity. These preference-driven models incorporate factors such as willingness to pay, travel schedules, and behavioral responses, often embedding them in market-based or game-theoretic frameworks.

Recent studies have begun to incorporate user-side behavior \citep{shahriar2020,kim2023} through game-theoretic or auction-based models \citep{cao2018,kim2022}, where users respond to incentives or participate in energy markets.
For example, \cite{kim2022} proposed a distributed auction-based mechanism for energy trading between EVs and mobile charging stations, where EVs submit bids and MCSs act as auctioneers. While the study highlights the benefit of user-side incentives, it assumes a single-round interaction with fixed pricing and does not address time-slot scheduling or bilateral negotiation. In contrast, our work introduces a multi-stage framework that jointly optimizes pricing and scheduling under infrastructure constraints, enabling adaptive coordination between users and a centralized operator.

Complementing these market-based approaches, several studies have explicitly modeled user preferences and behavioral heterogeneity in EV charging decisions \citep{wang2021,zhang2023,sica2025}. 
For example, \cite{zhang2023} incorporated users’ distance-based preferences into a bi-level planning model to analyze their impacts on station utilization and congestion.
\cite{sica2025} used discrete choice experiments to estimate users’ preferences for charging at home, work, or public locations, highlighting the effects of cost, waiting time, and convenience.
\cite{habbal2024} developed a multi-attribute recommendation model that ranks charging stations based on users’ priorities on waiting time, cost, and nearby facilities.
\cite{wang2021} applied binary logit and latent class models to capture heterogeneity in user choices, identifying service-sensitive and pragmatic user segments shaped by risk attitudes and satisfaction.
These studies emphasize that user preferences, heterogeneity, and behavioral traits can substantially influence system outcomes, yet they are rarely internalized into centralized decision-making frameworks.

\subsection*{Two-sided market and user-driven charging mechanisms}

Beyond operator- and preference-based formulations, a smaller set of studies has explored user-driven mechanisms that enable more direct interaction between users and operators. These include two-sided markets, auction models, and negotiation frameworks where both parties actively participate in determining charging outcomes.

Recent studies have explored market-based mechanisms to coordinate EV charging through bilateral interactions between users and service providers. For example, \cite{gao2021} propose a price-based iterative double auction for charger sharing markets, where private charger owners and EV users form a two-sided market. Their mechanism allows EV users to submit bids and compute welfare-maximizing schedules while accounting for time and location constraints. However, the platform centrally determines both user–station matching and charging schedules, which differs from settings where users retain the final choice of which charging station to use.
While these works explicitly analyze strategic interactions and equilibrium outcomes using game-theoretic tools (e.g., \cite{zavvos2021comprehensive,hajibabai2022game,ge2024distributed}), our framework adopts a centralized optimization framework with utility-based user participation, rather than conducting equilibrium analysis.
Another line of research investigates bilateral negotiations in the context of grid–EV coordination. For instance, \cite{zhang2020} develop a two-stage negotiation model between EVs and the power grid for discharging price determination, incorporating EV travel behavior and agent-based learning. While this model demonstrates the value of iterative price negotiation, it focuses on discharging incentives and does not consider time-slot scheduling or centralized system-level optimization.
In addition to auction-based or bilateral market mechanisms, some studies have explicitly investigated price negotiation between EV users and operators. For example, \cite{wang2018ev} propose a multi-round price negotiation framework where users and station operators iteratively adjust their bids until a mutually acceptable charging price is reached. Their results show that such negotiation can improve transaction volume and user satisfaction while reducing costs. However, the model mainly focuses on bilateral bargaining at the station–user level without considering system-level scheduling, time-slot allocation, or infrastructure constraints.

\subsection*{Summary and research gap}

Research on EV charging coordination can be broadly grouped into three categories. First, operator-driven approaches optimize prices and schedules from a centralized perspective, ranging from time-of-use and real-time tariffs \citep{wu2020,wang2019,adetunji2023} to more advanced methods such as bi-level formulations for station operations \citep{cedillo2022}, infrastructure-constrained scheduling \citep{wu2023}, and reinforcement learning–based dynamic pricing \citep{zhao2021,zulfiqar2024}. These methods improve efficiency but generally assume homogeneous and passive user behavior, overlooking heterogeneity and acceptance. Second, market-based and preference-driven models explicitly incorporate user heterogeneity, embedding factors such as distance, waiting time, risk attitudes, and satisfaction into choice models, auctions, or planning frameworks \citep{shahriar2020,kim2023,sica2025,habbal2024}. While these studies highlight user diversity, they typically optimize only user-side utility and lack mechanisms to align with operator objectives or centralized coordination. Third, two-sided and user-driven mechanisms emphasize direct interaction between users and operators. Examples include double auctions where users submit bids and operators compute allocations \citep{gao2021}, or negotiation frameworks allowing iterative price adjustment between EVs and stations \citep{wang2018ev}. While such works demonstrate the potential of bidding and negotiation, they usually remain confined to bilateral bargaining at the station level, without mechanisms for coordinated scheduling of multiple users.

Taken together, these streams reveal a clear research gap. Despite increasing attention to user-side behavior, existing approaches lack a framework that integrates heterogeneous user preferences with operator objectives under a coordinated scheduling environment. Attempts at negotiation remain limited in scope and do not extend to multi-user allocation and scheduling. Addressing this gap, the present study develops a user-driven framework that models the negotiation process through structured user bids, operator counteroffers, and acceptance decisions within a centralized optimization model. This design enables coordinated multi-user scheduling while embedding user-side heterogeneity, thereby balancing user utility with operator efficiency.

\section{Methodology} \label{Section3}

\begin{figure}
  \centering
  \includegraphics[width=.7\textwidth]{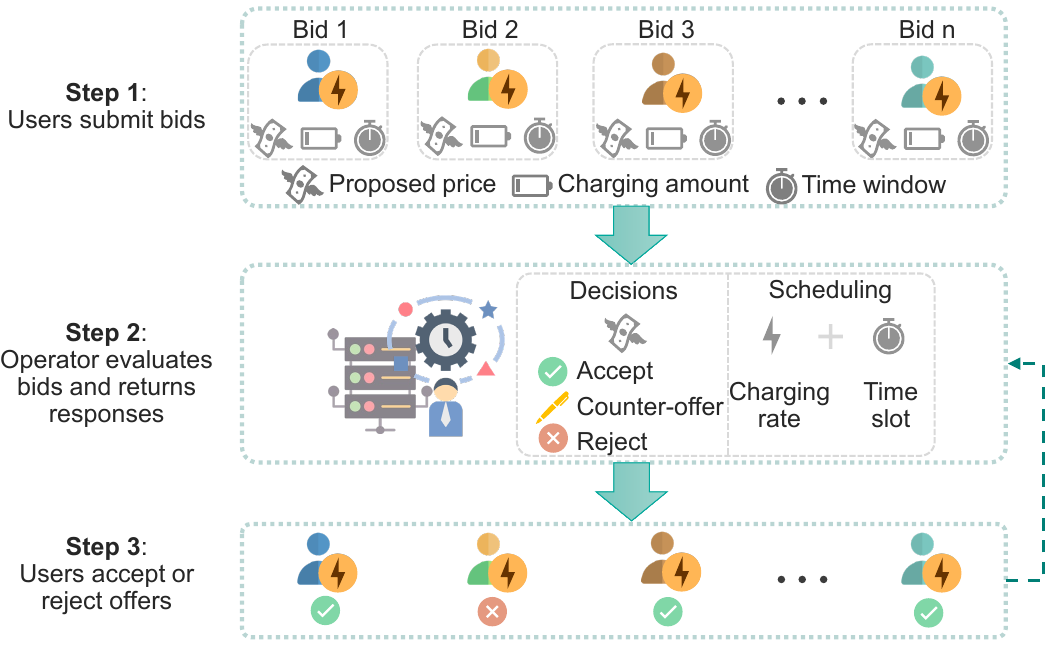}
  \caption{Three-step negotiation process between users and the operator.}\label{fig:negotiation-process}
\end{figure}

In this study, we propose a user-driven, bilateral negotiation framework with three sequential steps: user bidding, operator decision optimization, and user response, as illustrated in Figure \ref{fig:negotiation-process}. Specifically, users first submit charging bids that include their desired energy quantity, acceptable time windows, and the prices they are willing to pay. Based on the received offers and system constraints, the operator responds to offers with the aim of maximizing expected profits. Finally, EV users evaluate the received offers and decide whether to accept or reject them based on their utilities.

We consider a single public charging station operated over a finite planning horizon. The system is defined under the following assumptions:

\begin{itemize}
    \item Users submit their bids independently and do not have access to other users' bids. 
    \item The planning time horizon is divided into discrete time slots $t \in \mathcal{T}$ from which users specify their acceptable charging windows. 
    \item Different charging rates are available $\mathcal{R}$ (e.g., 22 kW, 50 kW, and 100 kW) at the charging station, and the operator decides how to allocate the chargers across users, considering compatibility.
    \item The implementation could be managed through a user-facing interface (e.g., a mobile app or web portal) that supports bidirectional interaction between users and the operator. The design of such interfaces is beyond the scope of this study.
\end{itemize}

This interactive process enables moderated pricing and adaptive scheduling that better aligns with both user preferences and system-level objectives. In the following sections, we elaborate on the three bidding and optimization steps, respectively. 

\subsection{Step 1: User Bidding}

Before the next planning horizon (e.g., the next day), the operator opens a bidding window, when each interested user $i$ submits a bid $s_i = (p_i^{\rm B}, Q_i^{\min}, Q_i^{\max}, \mathcal{T}_i)$. Each bid consists of three key elements characterizing the charging demand:

\begin{itemize}
    \item \textit{Proposed price $p_i^{\rm B}$:} Users enter the price they are willing to bid, per unit of energy. This reflects their urgency or preference for charging.
    
    \item \textit{Charging quantity in the range of $[Q_i^{\min}, Q_i^{\max}]$:} Users specify the minimum and maximum acceptable energy they would like to receive. In practice, this may be pre-filled or selectable based on their vehicle’s battery state and preferences (e.g., sliders or presets such as "Minimum 20 kWh, Maximum 40 kWh").
    
    \item \textit{Acceptable time window $\mathcal{T}_i$:} Users select their available time window within the discretized time horizon $t \in \mathcal{T}$, indicating when they would like to charge their EVs.
\end{itemize}

Given limited charging capacity and overlapping preferred time windows, users may perceive competitive pressure when formulating their bids. In this study, bids are treated as exogenous inputs that reflect users’ willingness to pay and time flexibility, and strategic equilibrium behavior among users is not explicitly modeled.
Each user naturally aims to strike a balance between offering a sufficiently attractive bid and satisfying their own flexibility constraints. Whether a bid is ultimately accepted depends not only on the competing bids submitted by other users, $\{s_j\}_{j \neq i}$, but also on the operator’s centralized allocation strategy, which is detailed in the following subsection.

\subsection{Step 2: Operator Decision Model for Bid Response}

After collecting all submitted bids from users for the upcoming time horizon, the operator performs centralized optimization to determine slot-level charging schedules and pricing decisions under operational constraints. The aim is to maximize profits while satisfying operational constraints (e.g., charger capacity and time slot availability) and ensuring rational pricing. Specifically, the operator determines whether to accept each user's bid, propose a counter-offer with a new price, or reject the request. For each user whose bid is accepted or countered, the operator further determines the charging rate $r \in \mathcal{R}$ and the allocated continuous time slots $t \in \mathcal{T}_i$.

This decision-making process is formulated as a mixed-integer linear program that jointly optimizes scheduling and pricing decisions as follows. Table~\ref{parameters} in the Appendix provides a detailed definition of each variable.

\begin{equation}
\max \sum_{i,t,r} (p^{\text{F}}_{i,t,r} - c_{r,t}) \cdot E_r \cdot z_{i,t,r}
\label{eq:obj}
\end{equation}

\noindent\textit{Subject to:}

\begin{equation}
\sum_r y_{i,r} = a_i \quad \forall i \in \mathcal{I}
\label{eq:acceptorreject}
\end{equation}

\begin{equation}
\sum_r z_{i,t,r} \leq A_{i,t} \cdot a_i \quad \forall i \in \mathcal{I}, t \in \mathcal{T}_i
\label{eq:availability}
\end{equation}

\begin{equation}
\sum_{t} z^s_{i,t} \leq 1,\quad 
  \sum_{t} z^e_{i,t} \leq 1 \quad \forall i,r
\label{eq:start-end-slot}
\end{equation}

\begin{equation}
z^s_{i,t} \geq z_{i,t,r} - z_{i,t-1,r} \quad \forall i,t,r
\label{eq:start-slot}
\end{equation}

\begin{equation}
z^e_{i,t} \geq z_{i,t,r} - z_{i,t+1,r} \quad \forall i,t,r 
\label{eq:end-slot}
\end{equation}

\begin{equation}
\sum_t q_{i,t} \geq Q^{\min}_i \cdot a_i, \quad \sum_t q_{i,t} \leq Q^{\max}_i \cdot a_i \quad \forall i
\label{eq:quantity-range}
\end{equation}

\begin{equation}
\sum_i q_{i,t} \leq C_t \quad \forall t
\label{eq:capacity}
\end{equation}

\begin{equation}
\sum_i z_{i,t,r} \leq M_r \quad \forall t,r
\label{eq:slot-limit}
\end{equation}

\begin{equation}
q_{i,t} = \sum_r E_r \cdot z_{i,t,r} \quad \forall i,t
\label{eq:quantity-from-slot}
\end{equation}

\begin{equation}
\sum_r z_{i,t,r} \leq y_{i,r}, \quad \sum_r y_{i,r} = 1 \quad \forall i,t
\label{eq:region-selection}
\end{equation}

\begin{equation}
x^{\text{B}}_{i,t,r} + x^{\text{N}}_{i,t,r} = z_{i,t,r} \quad \forall i,t,r
\label{eq:split-bn}
\end{equation}

\begin{equation}
p^{\text{N}}_{i,t,r} \leq (p^{\text{B}}_i + \delta_r) \cdot x^{\text{N}}_{i,t,r} \quad \forall i,t,r
\label{eq:new-price-cap}
\end{equation}

\begin{equation}
p^{\text{F}}_{i,t,r} = x^{\text{B}}_{i,t,r} p^{\text{B}}_i + x^{\text{N}}_{i,t,r} p^{\text{N}}_{i,t,r} \quad \forall i,t,r
\label{eq:final-price}
\end{equation}

\begin{equation}
p^{\text{F}}_{i,t,r} \geq (1+\varepsilon) \cdot c_{r,t} \cdot z_{i,t,r} \quad \forall i,t,r
\label{eq:min-profit}
\end{equation}

\begin{equation}
p^{\text{F}}_{i,t,r} \leq \alpha \cdot c_{r,t} \quad \forall i,t,r
\label{eq:max-price-multiplier}
\end{equation}

\begin{equation}
p^{\text{R}}_{t,r} \geq p^{\text{B}}_i \cdot x^{\text{B}}_{i,t,r} \quad \forall i,t,r
\label{eq:ref-price-lower-bound}
\end{equation}

\begin{equation}
p^{\text{N}}_{i,t,r} \geq p^{\text{R}}_{t,r} - p^{\text{H}}(1 - x^{\text{N}}_{i,t,r}) \quad \forall i,t,r
\label{eq:new-price-lower-than-ref}
\end{equation}

\begin{equation}
\sum_{i,t,r} x^{\text{B}}_{i,t,r} \geq \gamma \cdot \sum_{i,t,r} z_{i,t,r}
\label{eq:baseline-min-share}
\end{equation}

\begin{equation}
x^{\text{B}}_{i,t,r}, x^{\text{N}}_{i,t,r}, z_{i,t,r}, y_{i,r}, z^{\text{s}}_{i,t}, z^{\text{e}}_{i,t} \in \{0,1\}
\label{eq:binary-vars}
\end{equation}

\begin{equation}
p^{\text{N}}_{i,t,r}, p^{\text{F}}_{i,t,r}, q_{i,t} \geq 0
\label{eq:nonnegative-vars}
\end{equation}

Eq.~\eqref{eq:obj} maximizes the operator’s total profit, calculated as the sum of price-cost margins $p^{\text{F}}_{i,t,r} - c_{r,t}$ over all assigned charging slots $(i,t,r)$, weighted by the energy delivered $E_r$.
Eq. \eqref{eq:acceptorreject} ensures that if user $i$ is accepted to be served (either according to the original bid or after negotiation), then a charging rate is assigned to it.
Eqs. \eqref{eq:availability}–\eqref{eq:end-slot} are slot feasibility and charging continuity constraints. The user’s acceptable time window $\mathcal{T}_i \subseteq \mathcal{T}$ is encoded using a binary indicator $A_{i,t}$, where $A_{i,t} = 1$ if $t \in \mathcal{T}_i$ and $0$ otherwise. Eq. \eqref{eq:availability} ensures that users can only be assigned to time slots they declared acceptable. Eqs.~\eqref{eq:start-end-slot}–\eqref{eq:end-slot} model the continuity of charging by introducing binary variables to mark the start and end of charging, ensuring that each user has at most one continuous charging window. 
Specifically, Eq.~\eqref{eq:start-end-slot} limits each user to a single start and end slot at a given rate. Eqs.~\eqref{eq:start-slot}–\eqref{eq:end-slot} detect the start and end of charging by comparing adjacent time slots, thereby enforcing a single continuous charging window for each user. For boundary handling, $z_{i,0,r}$ and $z_{i,T+1,r}$ are defined as zero.
Eqs. \eqref{eq:quantity-range}–\eqref{eq:quantity-from-slot} are energy demand and capacity constraints. Eq. \eqref{eq:quantity-range} enforces that the total energy delivered to each user lies within their individual demand bounds. Eq. \eqref{eq:capacity} ensures that the total energy supplied in each slot does not exceed the system-wide power capacity. Eq. \eqref{eq:slot-limit} limits the number of users simultaneously using chargers at each rate. Eq. \eqref{eq:quantity-from-slot} defines the actual energy delivered to each user in each slot as a function of the assigned rate.
Eqs. \eqref{eq:region-selection}–\eqref{eq:split-bn} enforce consistency between charging-rate selection and slot-level assignment. Each user is assigned exactly one charging rate (Eq. \eqref{eq:region-selection}), and each assigned slot is classified as either accepted at the user’s original bid or assigned a counterprice (Eq. \eqref{eq:split-bn}).

Eqs. \eqref{eq:new-price-cap}–\eqref{eq:baseline-min-share} enforce pricing and rationality constraints. Eq. \eqref{eq:new-price-cap} limits the operator’s new price to a bounded markup above the original bid. Eq. \eqref{eq:final-price} defines the final price as either the accepted bid or a newly proposed price. Eq. \eqref{eq:min-profit} requires that any accepted slot generates a minimum profit margin. Eq. \eqref{eq:max-price-multiplier} caps the final price by a maximum allowed markup over cost. Eq. \eqref{eq:ref-price-lower-bound} sets the reference price for each slot and rate as the highest accepted bid. Eq. \eqref{eq:new-price-lower-than-ref} ensures that the newly proposed price is no lower than the reference price 
through a standard big-M formulation, where $p^{\text{H}}$ is a sufficiently large upper bound on feasible charging prices.
Eq. \eqref{eq:baseline-min-share} ensures that at least a $\gamma$ proportion of all assigned slot–rate pairs are priced at the user’s original bid price, thereby bounding the operator’s markup aggressiveness. Here, $\gamma$ denotes the price protection ratio, enforcing that no less than $\gamma$ fraction of accepted bids are served at their original prices. A larger $\gamma$ reflects stricter price preservation and reduced pricing flexibility.
Finally, Eqs. \eqref{eq:binary-vars}–\eqref{eq:nonnegative-vars} define variable domains, specifying binary decision variables and non-negativity conditions for all price and energy variables.

\subsection{Step 3: User Response and Acceptance Decision}

After receiving the operator’s offer, each user evaluates whether to accept the assigned charging plan, which includes the final price \( p_{i,t}^{\text{F}} \), the allocated charging rate \( r \in \mathcal{R} \), and the assigned time window \( \mathcal{T}_i^{\text{A}} \subseteq \mathcal{T}_i \). The user computes the utility from the offer and compares it with the utility of her best alternative charging arrangement (the outside option), such as charging at another public station or at home. Specifically, user \( i \) accepts the offer if and only if:
\begin{equation}
U_i(p_{i,t}^{\text{F}}, Q_{i,t}, r, \mathcal{T}_i^{\text{A}}) \geq U_i^{\text{O}},
\label{eq:utility}
\end{equation}
where $U_i^{\text{O}}$ denotes the utility associated with user $i$’s outside option. Otherwise, the user opts out, and no transaction occurs.

In this study, we adopt a static, utility-based acceptance rule: user $i$ accepts the operator’s offer if the realized utility exceeds their outside option $U_i^{\text{O}}$. The decision is myopic and no dynamic learning or intertemporal strategic behavior is modeled. For numerical illustration, we specify a representative utility function in which utility increases with the allocated charging rate and decreases with the markup over the submitted bid.
Heterogeneous preference parameters are introduced to capture user diversity in service valuation and price sensitivity.

This final step concludes the bilateral mechanism, where the user's participation is determined endogenously based on the trade-off between offered value and external alternatives.

\section{Numerical experiments} \label{Section4}

This section presents a series of numerical experiments designed to evaluate the proposed pricing and allocation model across diverse user and system conditions. The experiments are grouped by user characteristics (homogeneous vs. heterogeneous) and charger configurations (single vs. multiple charging rates), with an additional large-scale case to assess model scalability. All instances are solved to optimality within seconds using a commercial solver (i.e., Gurobi), reflecting the model’s computational tractability and structured formulation. These results underscore the framework’s potential for real-time integration in charging systems without relying on complex algorithmic interventions.

\subsection{Scenario design}

Representative scenarios are constructed to evaluate the model’s performance in responsiveness and fairness. Specifically, we vary the parametric dimensions, including user bid ranges, charger rate configurations, and pricing constraints. Key model parameters and baseline scenario settings are described below.

In the baseline scenario, user bid prices are drawn from a uniform distribution within 2.0–4.0 SEK/kWh, representing reasonable values below the prevailing public charging prices in Sweden, typically 4.0–8.0 SEK/kWh \citep{MerPrices2025, FortumIonity2024}. In extended scenarios, wider or biased bid distr ibutions (e.g., 1.0–6.0 or user-type-specific ranges) are used to examine the impacts of more heterogeneous user behaviors.

The model considers three prevailing charging rates, \( R = [22, 50, 100] \) kW, corresponding to slow, medium, and fast charging levels. Each time slot is 15 minutes long, and the unit cost per charging slot $c_{r,t}$ is fixed for each rate, with values of 1.5, 2.0, and 2.5 SEK for the three charging levels, respectively.

The maximum allowed markup over user bids is defined as \( \delta_r = [2, 2, 2] \) SEK, meaning that the operator is not allowed to increase the final price by more than 2 SEK per slot, regardless of the charging rate. This restriction is imposed to prevent potential abuse by the operator. Additionally, the system imposes a global maximum price-to-cost ratio of \( \alpha = 2.5 \) and enforces a minimum profit margin \( \varepsilon = 0.5 \) to ensure basic economic viability. These parameter values are chosen to reflect close-to-reality conditions in Sweden. However, the proposed model itself is general and can be applied to other contexts. 

\subsection{Homogeneous users}

We first examine the system behavior under homogeneous user profiles as a baseline scenario. In this setting, all users share identical charging demands ($Q^{\min}$, $Q^{\max}$), availability windows ($\mathcal{T}_i$), and submit randomly drawn bid prices within a predefined price range. This simplification allows us to attribute observed system behavior to the pricing mechanisms without interference from user-specific characteristics. It serves as a benchmark for understanding how operator strategies interact with user bids under idealized conditions.

\subsubsection{Single charging rate}

The operator adopts a single-rate pricing policy, meaning all users are assigned the same charging speed if accepted. Under this policy, the price protection ratio $\gamma$ is varied across scenarios to evaluate how strict bid preservation constraints affect overall system outcomes. This scenario provides a controlled environment to assess the fundamental trade-offs between price fairness, operator revenue, and user acceptance before introducing additional complexities such as heterogeneous charger types or user availability.

\noindent \textit{\textbf{Aggregated performance under random bidding}}

To examine the robustness of system outcomes under heterogeneous bidding realizations, we evaluate each $(\gamma, r)$ configuration over 50 independently generated scenarios, where bids are drawn randomly from the same distribution. This aggregation allows us to characterize
system-level trends that are independent of any single realization, providing statistically meaningful comparisons across charging rates and values of the price protection parameter $\gamma$. For each configuration, we report the mean and standard deviation of three key performance indicators: (i) operator profit, (ii) user acceptance rate, and (iii) the average final price paid by accepted users. In addition, we summarize pricing behavior using the average markup (defined as the difference between final and submitted bid prices for accepted users) and assess fairness based on the Gini coefficient of final unit prices.

Table~\ref{tab:summary-aggregate} presents the aggregated results from the 50 randomly generated bidding scenarios.
The Gini coefficient is computed using the standard formulation based on the dispersion of final unit prices paid by users whose requests are accepted by the system, and is then averaged over the 50 randomly generated bidding scenarios.
Several consistent patterns emerge. The user-driven framework behaves distinctly under low-, medium-, and high-power charging rates.
At 22 kW, both acceptance and final prices remain stable across $\gamma$ levels, and the resulting allocations exhibit high fairness.
At 50 kW, the system becomes more sensitive to $\gamma$, with acceptance rates declining and markups showing greater variability as price-preservation constraints tighten.
At 100 kW, users systematically underbid relative to the operator’s cost, rendering most requests economically infeasible. Under high $\gamma$, the operator is effectively unable to apply markups, leading to very low acceptance and near-zero markup.
This behavior can be intuitively explained by the bounded pricing range imposed by the model.
Since the allowed markup is defined as a percentage relative to the operator’s cost, lower charging rates correspond to a narrower absolute price range, resulting in limited variability across different values of $\gamma$.
At higher charging rates, the same proportional bounds translate into a wider absolute price range, making the system more sensitive to $\gamma$ and leading to greater variability in acceptance and pricing outcomes.

\begin{table}[width=.9\linewidth,cols=7,pos=h]
\centering
\caption{System performance over 50 random bidding scenarios (values shown as mean~$\pm$~std).}
\label{tab:summary-aggregate}
\begin{tabular*}{\tblwidth}{@{} LLLLLLL@{} }
\hline
Rate (kW) & $\gamma$ & Profit & Acceptance & Final price & Markup & Gini \\
\hline
22 & 0.2 & 291.76$\pm$3.77 & 0.60$\pm$0 & 3.71$\pm$0.09 & 0.65$\pm$0.63 & 0.01 \\
22 & 0.4 & 280.50$\pm$8.25 & 0.6$\pm$0 & 3.62$\pm$0.17 & 0.49$\pm$0.62 & 0.02 \\
22 & 0.6 & 262.99$\pm$13.02 & 0.6$\pm$0 & 3.49$\pm$0.27 & 0.37$\pm$0.60 & 0.04 \\
22 & 0.8 & 239.82$\pm$17.42 & 0.6$\pm$0 & 3.32$\pm$0.36 & 0.19$\pm$0.50 & 0.05 \\
\hline
50 & 0.2 & 620.49$\pm$26.30 & 1.00$\pm$0 & 4.48$\pm$0.50 & 1.51$\pm$0.78 & 0.06 \\
50 & 0.4 & 511.56$\pm$74.45 & 0.94$\pm$0.13 & 4.19$\pm$0.60 & 1.17$\pm$0.97 & 0.08 \\
50 & 0.6 & 350.91$\pm$96.89 & 0.72$\pm$0.19 & 3.95$\pm$0.64 & 0.75$\pm$0.97 & 0.09 \\
50 & 0.8 & 238.72$\pm$77.72 & 0.55$\pm$0.18 & 3.74$\pm$0.55 & 0.35$\pm$0.73 & 0.07 \\
\hline
100 & 0.2 & 272.02$\pm$228.13 & 0.47$\pm$0.41 & 4.82$\pm$0.65 & 1.60$\pm$0.80 & 0.07 \\
100 & 0.4 & 127.69$\pm$125.06 & 0.23$\pm$0.22 & 4.76$\pm$0.81 & 1.13$\pm$1.00 & 0.09 \\
100 & 0.6 & 59.89$\pm$71.47 & 0.13$\pm$0.14 & 4.31$\pm$0.77 & 0.52$\pm$0.88 & 0.04 \\
100 & 0.8 & 33.45$\pm$30.43 & 0.10$\pm$0.09 & 3.87$\pm$0.07 & 0$\pm$0 & 0.00 \\
\hline
\end{tabular*}
\end{table}

\begin{figure}
  \centering
  \includegraphics[width=.7\textwidth]{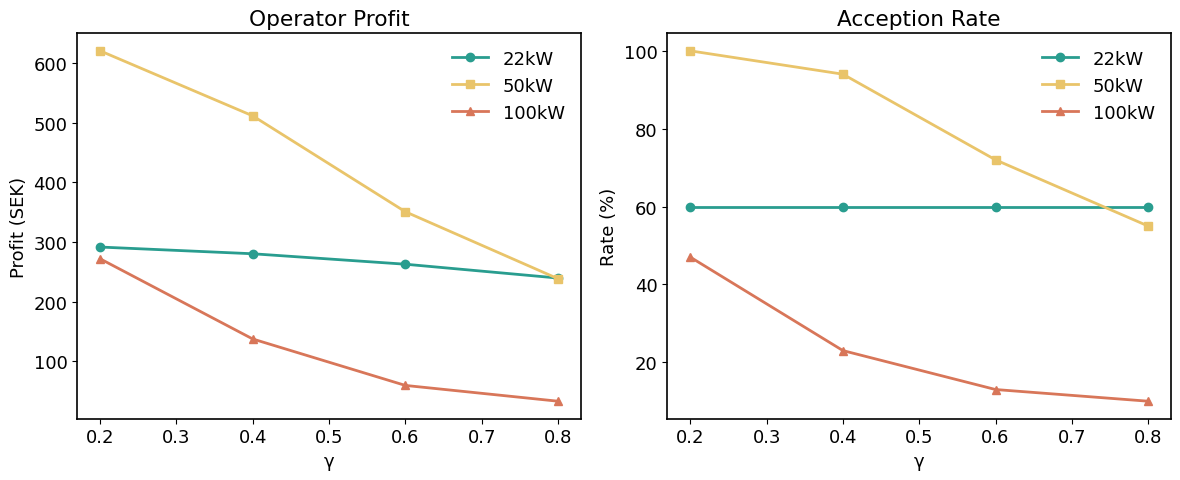}
  \caption{Impact of $\gamma$ on system outcomes.}\label{fig:system-outcome}
\end{figure}

To illustrate the most dominant system-level trends, Figure \ref{fig:system-outcome} visualizes the aggregated mean operator profit and acceptance rates across the price protection ratio $\gamma$ and charging rates. These two indicators highlight the main behavioral responses to changes in $\gamma$, while the remaining performance metrics are summarized in Table \ref{tab:summary-aggregate}.

\noindent \textit{\textbf{Representative scenario visualization}}

\begin{figure}
  \centering
  \begin{subfigure}[b]{.9\textwidth}
    \centering
    \includegraphics[width=\textwidth]{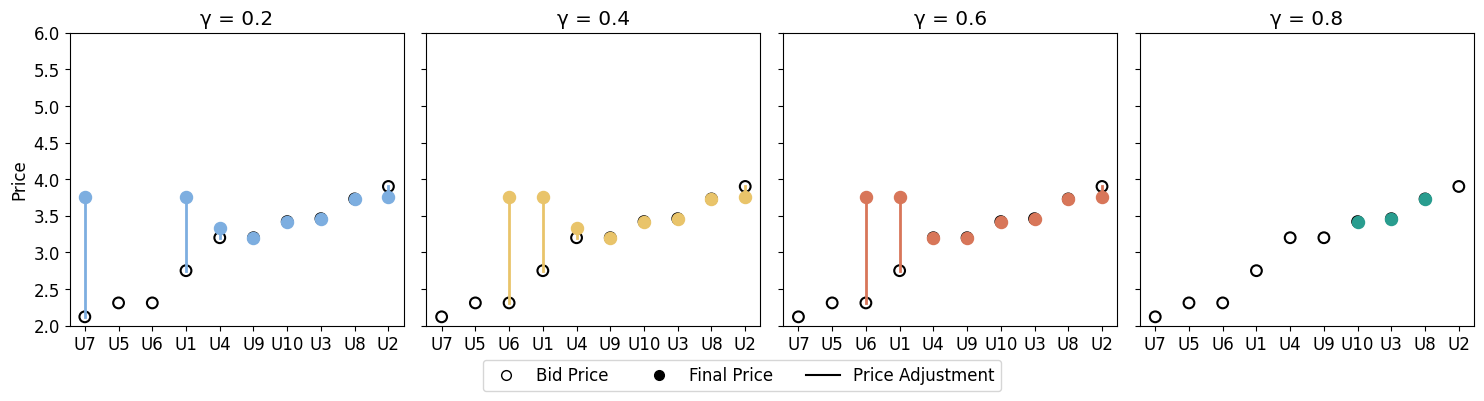}
    \caption{22kW}
  \end{subfigure}

  \begin{subfigure}[b]{.9\textwidth}
    \centering
    \includegraphics[width=\textwidth]{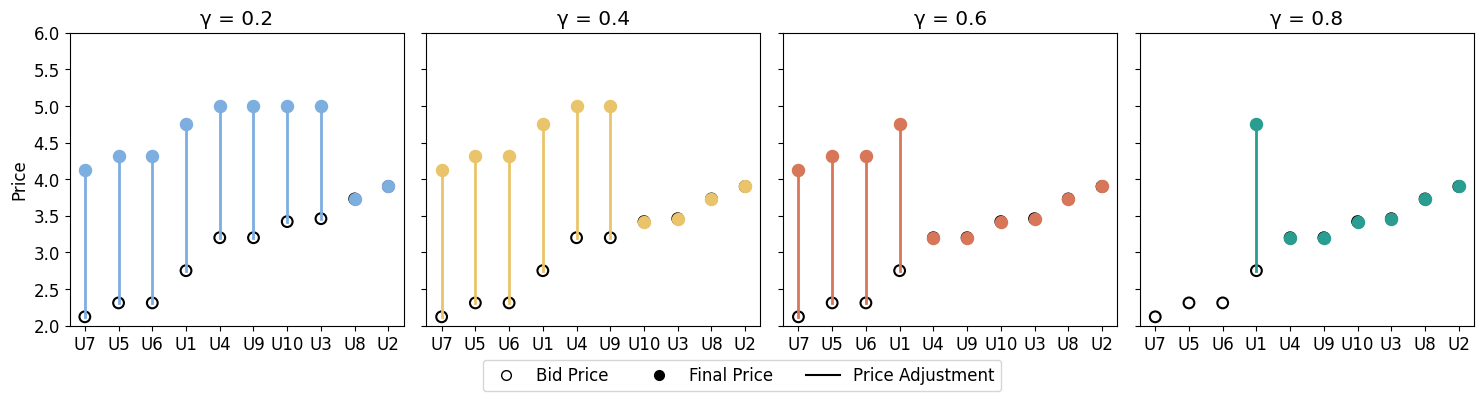}
    \caption{50kW}
  \end{subfigure}

  \begin{subfigure}[b]{.9\textwidth}
    \centering
    \includegraphics[width=\textwidth]{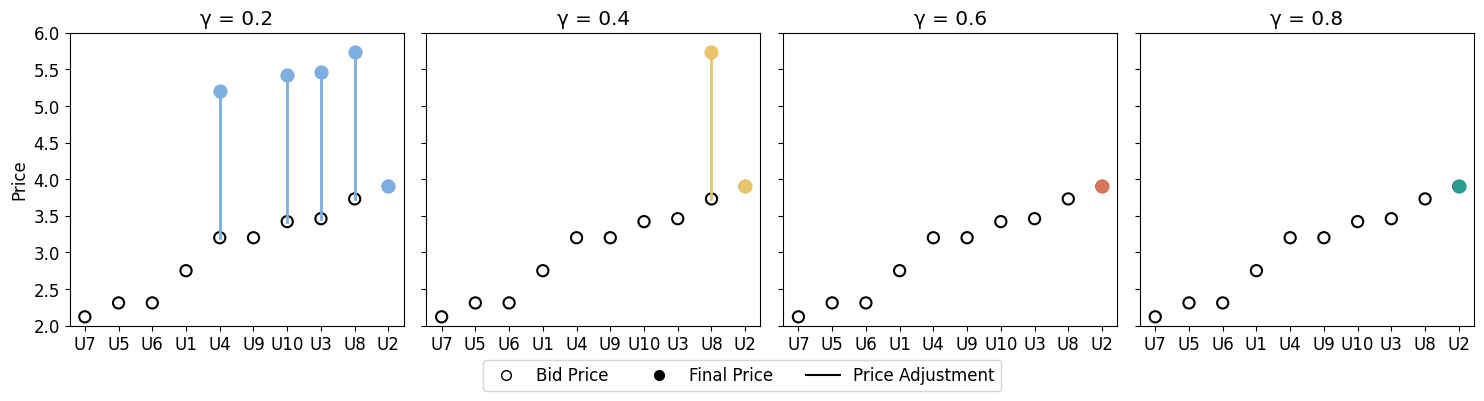}
    \caption{100kW}
  \end{subfigure}

  \caption{Final prices vs. user bidding prices under different $\gamma$ values.}
  \label{fig:final-vs-bid}
\end{figure}

To complement the aggregated statistics, we select representative scenarios to illustrate user-level pricing and acceptance outcomes under different price protection ratios, ~$\gamma$.
Figure~\ref{fig:final-vs-bid} shows the results for the single-rate cases. It compares users’ accepted prices with their submitted bids across different levels of price protection~($\gamma$), isolating the effect of pricing flexibility at each charging rate. This provides a clear view of how stricter bid preservation constrains price adjustments and influences user acceptance behavior.
In each panel, open circles denote users’ bid prices, while filled circles indicate the accepted final prices under the corresponding value of~$\gamma$. The vertical line segment between the two markers represents the markup. Users with only an open circle (i.e., no line or filled marker) are rejected.

As shown in Figure~\ref{fig:final-vs-bid}, filled circles represent users' final prices under a given $\gamma$, while open circles denote their bid prices. Vertical lines indicate the markup between the bid and the accepted price. The absence of a filled circle indicates that the corresponding user is rejected. For the 22 kW system, a fixed group of users is consistently rejected across all $\gamma$ values, reflecting capacity limitations. Among those served, accepted prices show significant markups at low $\gamma$, which progressively shrink as $\gamma$ increases. This suggests that $\gamma$ mainly regulates price markups, while service eligibility is governed by system capacity. In the 50 kW case, most users are served under low to moderate $\gamma$. Rejections emerge only at $\gamma=0.8$, indicating a threshold beyond which the system can no longer afford to accept high bids. Price adjustments remain moderate and clearly respond to $\gamma$, showcasing a balance between capacity and pricing flexibility. In contrast, the 100 kW system exhibits a sharp drop in accepted users as $\gamma$ increases. While full acceptance and aggressive markups are seen at $\gamma=0.2$, the system rapidly loses its ability to serve users once pricing is constrained. This illustrates a key vulnerability: high-power chargers rely heavily on pricing flexibility to remain profitable and inclusive. Overall, the results demonstrate how the interaction between charger capacity and pricing constraints jointly determines who gets served and at what price. Capacity limits define the upper bound of service coverage, while $\gamma$ defines the margin for revenue recovery via price adjustments.

From a user-centered perspective, the bid price has a clear influence on both the likelihood of acceptance and the final price paid. Users with higher bids are consistently more likely to be directly accepted, especially under tight capacity or low price protection ($\gamma$). In contrast, users who bid too low are either rejected outright or subject to substantial markups, particularly when the system has the flexibility to adjust prices. Interestingly, users with moderate bids are the most likely to be counter-offered, they fall within the operator’s profitable range but are still low enough to trigger price adjustment. This suggests a strategic trade-off: bidding too low increases rejection risk, while bidding too high guarantees acceptance but forgoes cost savings. To maximize acceptance while minimizing cost, users may benefit from slightly above-average bids, especially when system congestion is high or $\gamma$ is low. Additionally, targeting non-peak time slots or periods of lower demand may help low-bidding users avoid rejection due to competition. These dynamics highlight the importance of transparent feedback or system learning to guide user bidding behavior over time.

\subsubsection{Multiple charging rates}

We next examine a more realistic setting where the operator offers multiple charging options (e.g., 22 kW, 50 kW, and 100 kW), and users may be assigned to any available charger based on pricing and scheduling optimization. All users still share homogeneous demand and time availability, but charger heterogeneity introduces an additional layer of decision-making for the operator. This setting reflects practical deployment scenarios and allows us to assess how the pricing mechanism adapts to variable charging infrastructure. In particular, we focus on how bid preservation constraints affect the assignment of charger types and resulting price fairness.

\begin{figure}
  \centering
    \includegraphics[width=\textwidth]{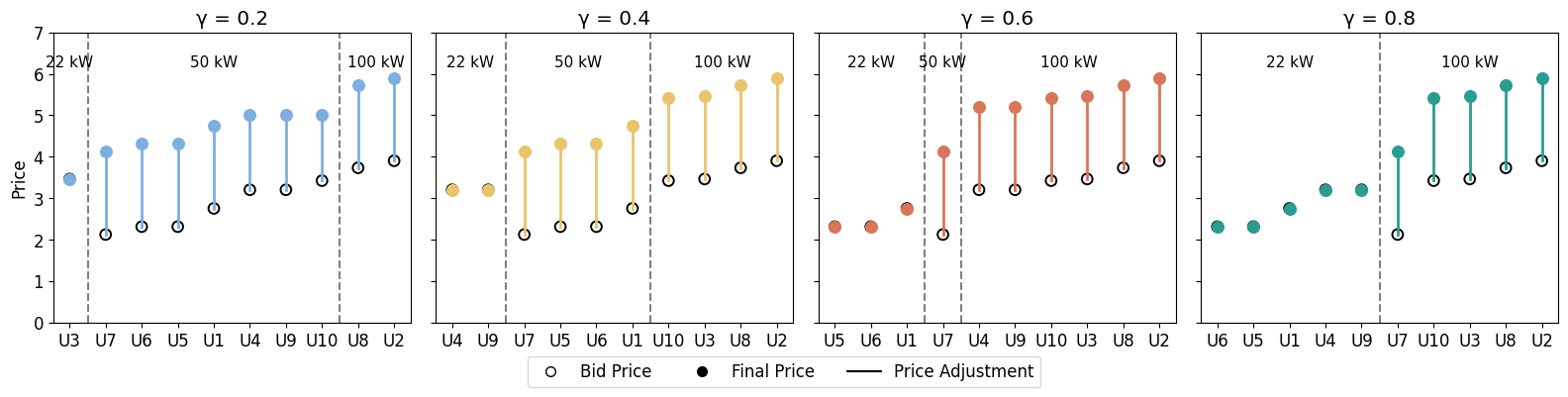}

  \caption{User-level outcomes under varying values of $\gamma$.}
  \label{fig:final-vs-bid-rate}
\end{figure}

Figure~\ref{fig:final-vs-bid-rate} presents user-level outcomes under varying values of $\gamma$. In each subplot, open circles denote users’ bids, while filled circles represent the corresponding accepted prices. Vertical lines indicate the markup between bids and final accepted prices. The dashed vertical separators distinguish users assigned to different charger capacities (22 kW, 50 kW, and 100 kW). The scenario assumes heterogeneous charger availability (one charger for each capacity level) and randomly drawn user bids. Across the four subplots ((a)--(d)), the bid preservation parameter $\gamma$ increases from 0.2 to 0.8, allowing an examination of its impact on pricing and rate allocation decisions.

At low $\gamma$ values (e.g., $\gamma = 0.2$), most users are charged a price above their original bid, reflecting the operator’s flexibility to optimize revenue through selective price increases. Higher charging rates (50\,kW and 100\,kW) are frequently assigned under this setting, as relaxed bid preservation allows the operator to prioritize profitability in the price-rate combination. As $\gamma$ increases, a clear shift occurs: final prices converge toward the bid values, while assignments of lower-rate chargers (22\,kW) become more common. By $\gamma = 0.8$, many users pay exactly their bid price but receive the lowest available rate. This pattern indicates that tighter bid preservation constraints reduce the operator’s pricing flexibility and lead to more conservative, resource-saving allocation strategies.

Importantly, the figure also reveals that being a relatively high bidder among users does not necessarily ensure preservation of the original bid or access to premium service. For example, \textit{user 2} and \textit{user 8} consistently offer relatively high bids, yet under all $\gamma$ settings, both are assigned 100\,kW chargers with final prices exceeding their bids. This suggests that their bids, while high compared to others, remain below the internal price threshold required to receive 100\,kW service without a counter-offer. In contrast, users with lower bids who are allocated 50\,kW or 22\,kW chargers often pay exactly their bid, particularly at higher $\gamma$ values. In addition, the mechanism promotes preference-aligned fairness: users who pay more do so in exchange for higher service levels, and price differences reflect self-selected urgency rather than arbitrary treatment.

Overall, these results highlight that the operator’s counter-offer decision is based on unrevealed \textit{absolute price thresholds} rather than relative user rankings. While bidding high increases the likelihood of receiving faster charging, it does not guarantee price preservation. As $\gamma$ increases, bid fairness improves and markups are suppressed---but often at the cost of reduced charger utilization efficiency. For users aiming to avoid counter-offers, understanding the internal pricing logic and service thresholds is more effective than simply outbidding peers.

\smallskip
\smallskip

\noindent\textit{\textbf{Expanded bid range}}

To further investigate how absolute bid levels affect pricing outcomes and service allocation, we extend the user bid range from the original interval \([2, 4]\) to a broader span of \([1, 6]\). This revised setting introduces greater diversity in bidding behavior, including aggressive underbidding and generous overbidding. The objective is to evaluate whether high absolute bids can effectively avoid counter-offers and secure premium charging rates, and whether low bids are more likely to result in rejections or downgraded service under varying bid preservation constraints.

\begin{figure}
  \centering
    \includegraphics[width=\textwidth]{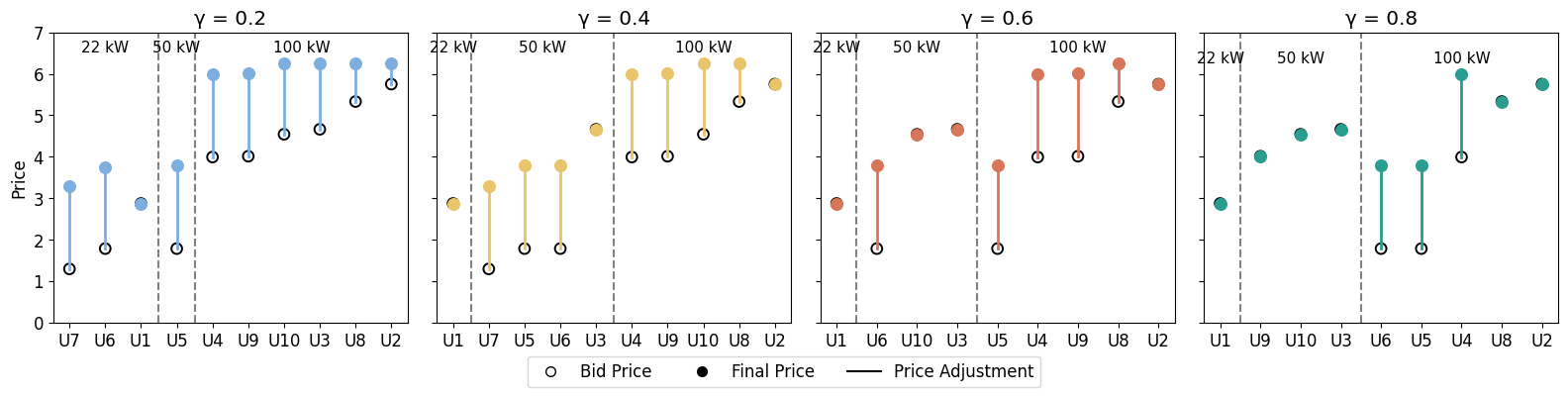}
  \caption{User-level outcomes under varying values of $\gamma$ with expanded bid range $[1, 6]$}
  \label{fig:expanded-final-vs-bid-rate}
\end{figure}

Figure~5 presents user-level outcomes under different values of $\gamma$, using the expanded bid range. In this setting, the role of absolute bid levels becomes more pronounced. High-bidding users, such as user 2 and user 4, consistently receive 100 kW charging across all $\gamma$ values with either no markup or only minor price increases, indicating their bids exceed the internal threshold for premium service. In contrast, user 1 regularly receives 22 kW charging without a price increase, suggesting the bid aligns closely with the internal threshold for low-rate service. Conversely, user 7 is rejected service under higher $\gamma$ values (e.g., $\gamma=0.6$ and $0.8$), indicating that their bid falls below the minimum price level acceptable for any service when constraints tighten. Overall, this experiment confirms that absolute bid levels, not relative rankings, govern both price preservation and service allocation.

\smallskip
\smallskip

\noindent\textit{\textbf{Different bid ranges}}

\begin{figure}
  \centering
    \includegraphics[width=\textwidth]{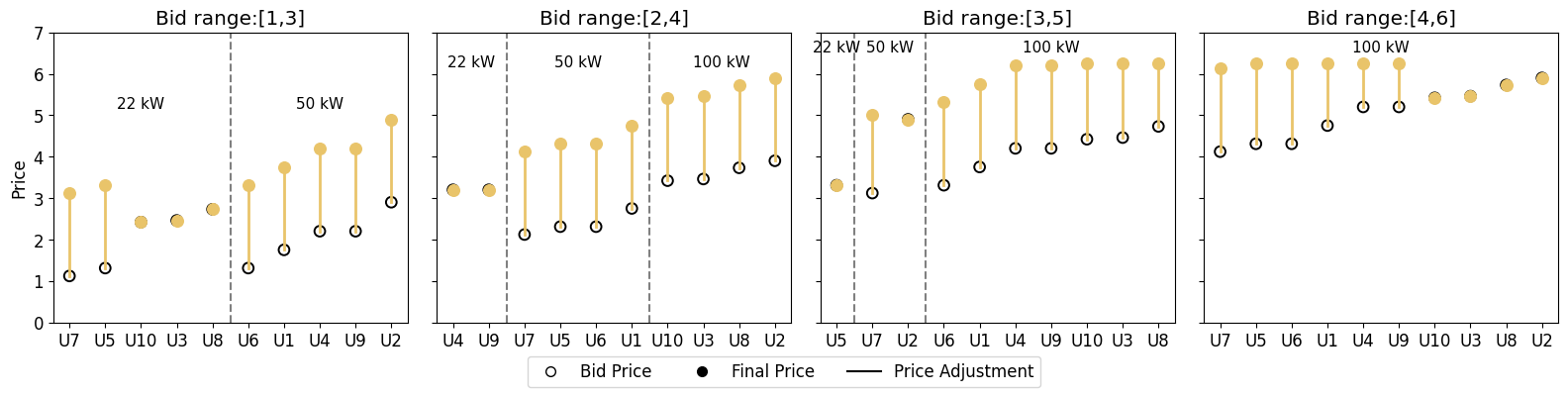}
  \caption{User-level outcomes with different bid range}
  \label{fig:differentbidrange}
\end{figure}

Following the expanded bid range experiment, we further investigate how the overall bidding environment affects user outcomes by fixing $\gamma=0.4$ and varying the bid range itself. The purpose is to isolate the effect of absolute bid magnitudes and examine whether a system-wide shift toward lower or higher bids changes price outcomes and charger assignments. A very low range such as $[0,2]$ produces no feasible solution under the operator’s internal thresholds and is therefore omitted. Figure \ref{fig:differentbidrange} shows user-level outcomes under several different bid ranges while keeping $\gamma=0.4$ fixed.

When users collectively bid lower, such as in the $[1,3]$ range, most bids fall below the thresholds required for premium service. As a result, users are still served, but predominantly with 22 kW chargers, and price adjustments remain small. This indicates that low bids do not necessarily lead to rejection; instead, they lead to downgraded service levels. As the bid range shifts upward to $[2,4]$ and $[3,5]$, more users cross the internal thresholds, leading to increased assignments of 50 kW and 100 kW chargers with moderate markups. In the highest range $[4,6]$, nearly all users qualify for premium service, resulting in uniform 100 kW allocations and minimal counter-offers. Overall, the results confirm that absolute bid levels shape service allocation: low collective bids still allow feasible service but mainly with lower-power chargers, whereas higher bid ranges consistently unlock premium charger assignments.

\subsection{Heterogeneous users}

To further examine how temporal heterogeneity affects pricing and allocation, we relax the assumption of identical availability windows and assign distinct time slot availabilities to each user. This introduces a new dimension of user heterogeneity, capturing more realistic charging behaviors in which users differ in both bid prices and scheduling flexibility. Specifically, we construct a scenario where the first five users are available throughout the entire scheduling horizon, while the remaining five users are only available during a narrow 2-hour window in the middle of the day. Both groups share the same bid structure, with user bids restricted to the interval [2, 4]. This setup allows us to isolate the effect of temporal flexibility from bidding behavior and observe how the scheduling mechanism prioritizes users under different levels of time availability.

\subsubsection{Single charging rate}

Similar to the previous section, we begin by analyzing the case with a single 22 kW charger. This allows us to observe how time-constrained users compete with flexible ones under different levels of pricing constraints.

\begin{figure}
  \centering
  \includegraphics[width=\textwidth]{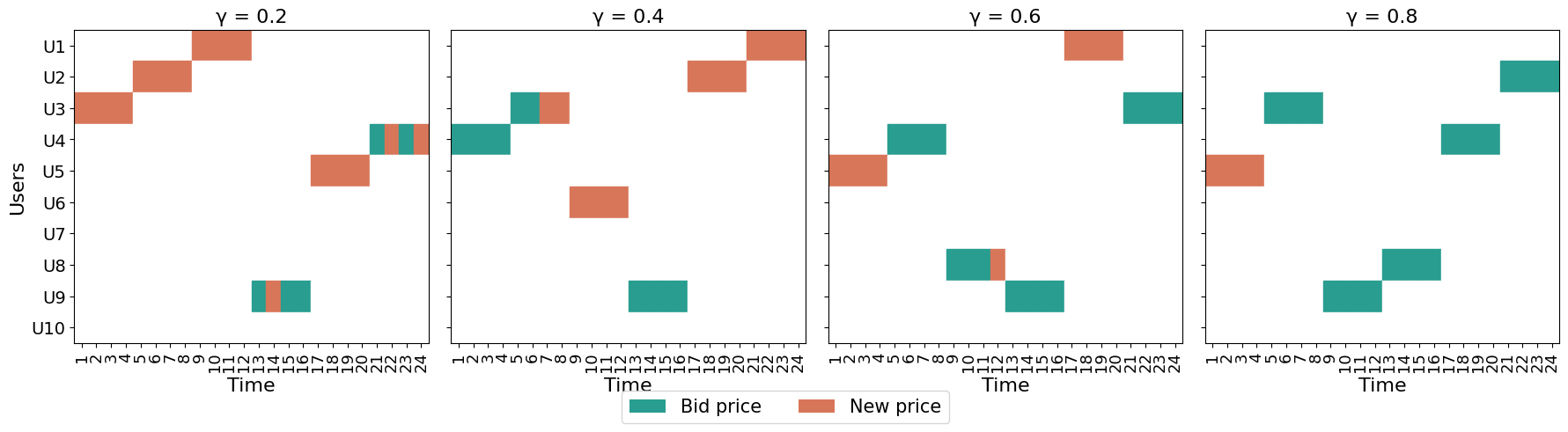}
  \caption{User-level allocation outcomes under under single 22 kW charger}\label{fig:allocation-outcome}
\end{figure}

Figure~\ref{fig:allocation-outcome} visualizes the user-slot allocation outcomes under varying bid preservation constraints \( \gamma \in \{0.2, 0.4, 0.6, 0.8\} \). A clear trade-off emerges: as \( \gamma \) increases, the system shifts from aggressive price adjustments (orange) to strict bid preservation (green), reducing pricing flexibility. When \( \gamma = 0.2 \), most allocations are made at higher prices, and time-flexible users (Users~1--5) are fully scheduled. As \( \gamma \) increases, these users begin to lose priority to high-bid but time-constrained users (e.g., Users~8 and~9), who can be served without price adjustments. At \( \gamma = 0.8 \), the system heavily favors bid-preserving allocations, excluding even flexible users if their bids are not sufficiently high. This demonstrates that while lower \( \gamma \) values allow revenue optimization through pricing, higher \( \gamma \) values improve price fairness at the cost of reduced scheduling flexibility and overall service coverage. The results highlight the critical interplay between temporal flexibility, bid level, and bid preservation constraints in determining allocation outcomes.

\begin{figure}
  \centering
  \includegraphics[width=\textwidth]{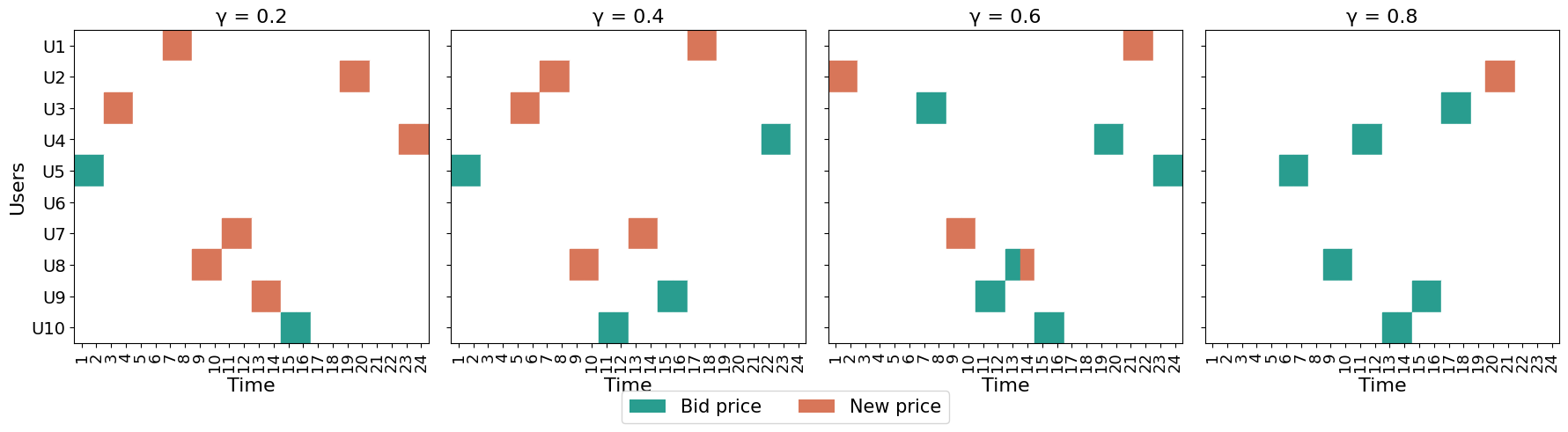}
  \caption{User-level allocation outcomes under under single 50 kW charger}\label{fig:allocation-outcome_50kW}
\end{figure}

Compared to the 22~kW case, the single 50~kW charger leads to a denser allocation structure, where each user typically occupies fewer time slots due to the higher power per slot. This compression in slot usage allows more users to be served within tight availability windows, improving scheduling feasibility under time constraints. However, since each slot represents a larger energy block, the system becomes more sensitive to bid levels: users with marginal bids may be excluded unless their flexibility aligns perfectly with the remaining capacity. As shown in Figure~\ref{fig:allocation-outcome_50kW}, while the single 50~kW charger enables denser scheduling due to its higher power per slot, increasing $\gamma$ values still impose notable bid preservation constraints. At $\gamma=0.8$, Users~1, 6, and 7, are entirely excluded from service, despite varying levels of temporal flexibility. This reflects a strict enforcement of bid thresholds, where even users with sufficient time availability may be denied access if their bids fall short. The result highlights that under strong bid preservation, the system shifts toward prioritizing high-bid users, such as Users~9 and 10, while sacrificing pricing flexibility and reducing user inclusiveness.

\subsubsection{Multiple charging rates}

\begin{figure}
  \centering
  \includegraphics[width=\textwidth]{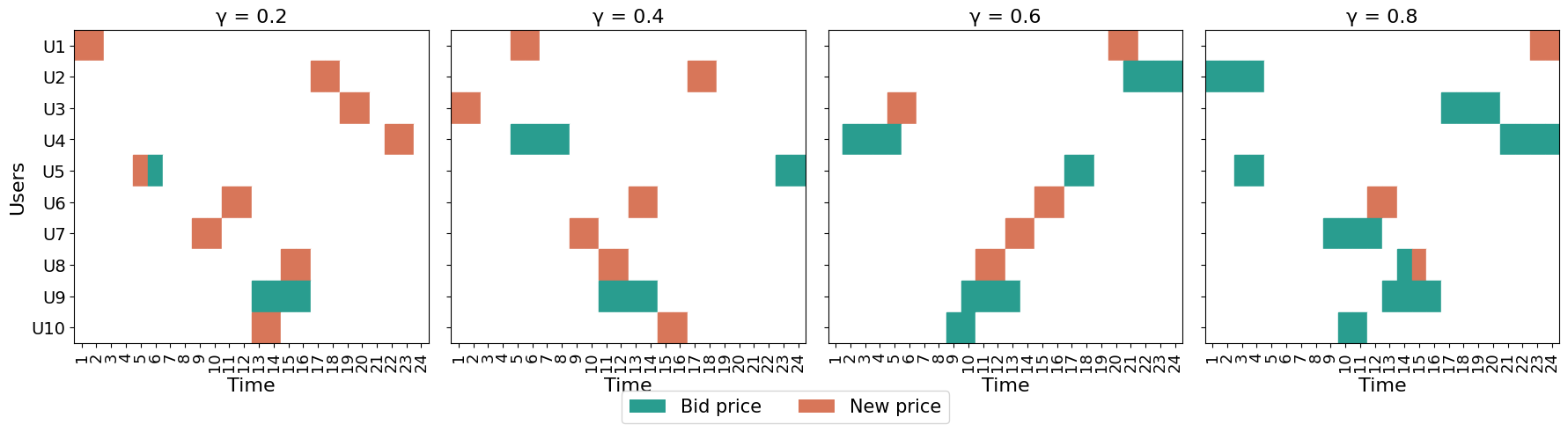}
  \caption{User-level allocation outcomes under one 22 kW charger and one 50 kW charger}\label{fig:allocation-outcome_2rates}
\end{figure}

We further consider a scenario where different charging rates are available. This reflects realistic infrastructure deployments in which both fast chargers (e.g., 50 kW) and regular chargers (e.g., 22 kW) coexist to serve a heterogeneous set of EVs. We examine how this multi-rate configuration influences allocation under different bid preservation levels. Figure \ref{fig:allocation-outcome_2rates} presents the user-level allocation outcomes under varying $\gamma$ values when one 22 kW charger and one 50 kW charger are jointly used to serve the system.

In the mixed-power scenario where one 22 kW charger and one 50 kW charger jointly serve the system, all users are successfully scheduled across all $\gamma$ levels. This setup includes both 22 kW and 50 kW users with heterogeneous bid prices and time availability. The outcome demonstrates that combining chargers of different rates enhances scheduling flexibility in a way that neither charger alone could achieve. Specifically, the 50 kW charger efficiently handles high-demand users with fewer time slots, while the 22 kW charger complements this by covering lower-power users or filling gaps in the schedule that might otherwise be unused due to compatibility or fairness constraints.
Compared to the single-charger scenarios, this configuration avoids slot congestion and reduces the need for aggressive price adjustments. Even under high bid preservation ($\gamma=0.8$), the system satisfies all users without violating constraints, indicating that heterogeneous charger configurations can improve both allocation efficiency and fairness. The joint use of different power levels also allows better matching between user demand and charger capability, enabling finer-grained resource utilization.
Overall, these results suggest that deploying a mix of charger types, rather than uniform infrastructure, can be a highly effective strategy in shared charging systems, particularly when balancing price constraints, time-window heterogeneity, and power demand diversity.

Compared to the single-charger scenarios, the joint 22 kW and 50 kW configuration offers clear improvements in both flexibility and inclusiveness. In the 22 kW-only setup (Figure~\ref{fig:allocation-outcome}), limited capacity prevents full user coverage under high $\gamma$ values. Similarly, while the 50 kW-only case (Figure~\ref{fig:allocation-outcome_50kW}) provides more power, it still excludes flexible users such as User 1 and User 7 when $\gamma=0.8$. In contrast, the mixed-rate configuration (Figure~\ref{fig:allocation-outcome_2rates}) consistently satisfies all users across all $\gamma$ values, highlighting how heterogeneous charger deployment mitigates both power bottlenecks and fairness constraints.

\subsection{Large-scale cases}

\begin{figure}
  \centering
 \includegraphics[width=.6\textwidth]{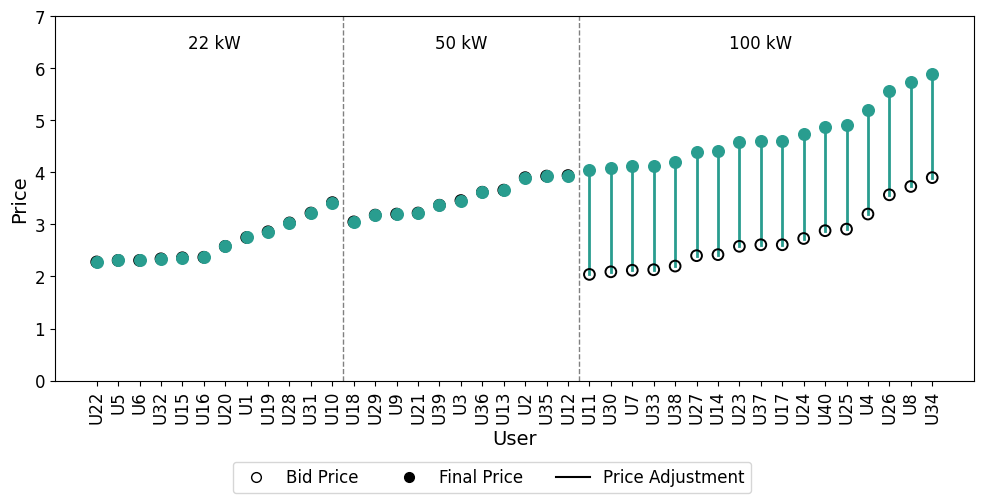}
  \caption{Bid and counter prices under a narrow bid range $[2,4]$ with $\gamma = 0.8$}\label{fig:40users_48slots_price}
\end{figure}

\begin{figure}
  \centering
 \includegraphics[width=.6\textwidth]{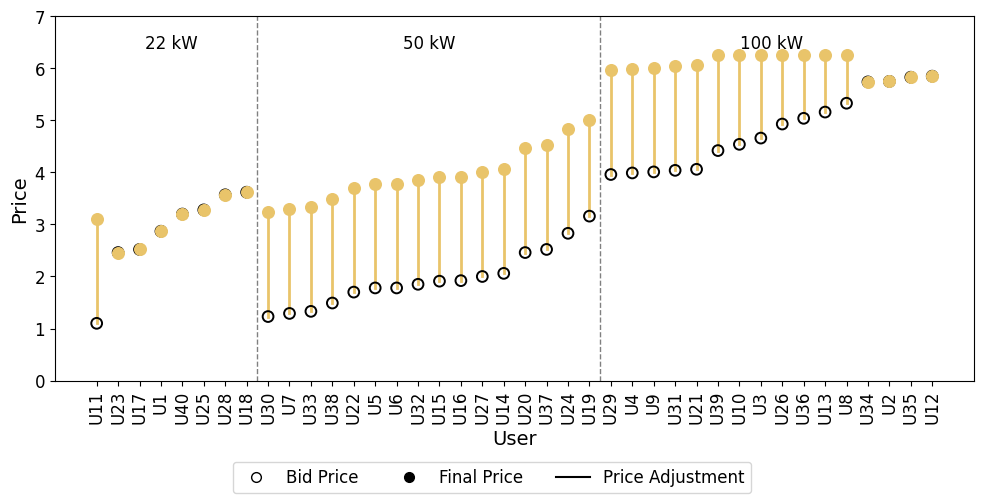}
  \caption{Bid and counter prices under a wide bid range $[1,6]$ with $\gamma = 0.4$}\label{fig:40users_48slots_price_wide}
\end{figure}

To gain a more comprehensive view of system-level behaviors, we further examine the full optimal schedule for 40 users across 48 time slots and three charger types (22 kW, 50 kW, and 100 kW). 
Instead of visualizing the full slot-level schedule, we focus directly on the resulting pricing patterns (Figures~\ref{fig:40users_48slots_price} and \ref{fig:40users_48slots_price_wide}), which reveal how bid ranges influence the uniformity of price adjustments and how charger heterogeneity determines the relative magnitude of counter-offers.

Figure~\ref{fig:40users_48slots_price} displays, for each user, the bid price (white circles) and the final price (green dots), sorted by charger type. Under the narrow bid range of $[2,4]$ and $\gamma=0.8$, all deviations between $p_u^{\mathrm{bid}}$ and $p_u^{\mathrm{final}}$ occur exclusively among users allocated 100~kW chargers. Out of 40 users, 17 (42.5\%) receive 100~kW service, and each experiences the same increment $\delta_r$ (equal to 2 in this setting). Bids of users served by 22~kW or 50~kW chargers are accepted without adjustment.
This pattern arises from the interaction between the fixed markup $\delta_r$ and the price cap. Within this bid range, none of the users approach the cap, so the markup can be applied without violating feasibility. However, applying the same markup to higher bids generates a larger absolute increase in revenue under the cap. Since users assigned to 100 kW chargers are located in the upper part of the bid distribution, they become the most profitable candidates for uniform markup application. In contrast, users served by 22 kW or 50 kW chargers yield lower incremental gains, making price adjustments less attractive. As a result, the narrow-range case exhibits a uniform uplift pattern that reflects profit prioritization rather than bid heterogeneity.

Figure~\ref{fig:40users_48slots_price_wide} displays the corresponding outcome under a wider bid range of $[1,6]$ and a moderate preservation ratio, $\gamma=0.4$. In this setting, the bid range becomes wide, and the interaction between the fixed markup and the price cap leads to price adjustments appearing across all charger types, unlike the narrow-range case. Users with the highest bids are already close to the feasible upper bound, so applying $\delta_r$ provides little additional revenue and risks violating the cap, resulting in no adjustment for these users. In contrast, users with low bids have substantial slack relative to the cap, making the markup feasible and profitable regardless of whether they are assigned 22~kW, 50~kW, or 100~kW service.
As a consequence, low-bid users across all charger types consistently receive the markup, whereas high-bid users retain their original bids. Mid-range users exhibit mixed outcomes: some receive the markup, while others do not, depending on whether adding $\delta_r$ remains feasible under the price-cap and the bid-preservation requirement~$\gamma$, as well as on whether they are allocated fast-charging service. This produces a heterogeneous adjustment pattern that contrasts sharply with the uniform uplift observed in the narrow-range scenario.

Overall, the two cases illustrate that counter-price patterns are not determined by bids alone, but emerge from the joint interaction of bid heterogeneity, the price-cap parameter~$\alpha$, and the design of the markup rule~$\delta_r$. The same fixed increment can generate either uniform adjustments or stratified outcomes, depending on how closely users' bids sit relative to the price cap and on how charging capacity is allocated.
These findings show that the markup mechanism is a structural component of the pricing model: by shaping which users ultimately receive counter-price adjustments, it governs how price flexibility is distributed across charger types and bid levels. This highlights the broader insight that counter-price design, through the choice of markup magnitude and its interaction with the price cap, plays a central role in balancing revenue, fairness, and capacity rationing in operator-driven charging markets.

\subsection{Utility-Based User Response Simulation}

To evaluate how users may respond to the operator’s counter-offers under heterogeneous preferences, we simulate utility-based acceptance behavior using the formulation introduced in Step 3. Specifically, each user \( i \) accepts the assigned charging plan if their utility defined as follows,
\begin{equation}
    U_i = \theta_{1,i} \cdot R_r - \theta_{2,i} \cdot \sum_{t \in \mathcal{T}_i^{\text{alloc}}} (p^{\text{F}}_{i,t} - p^{\text{B}}_{i,t}) Q_{i,t},
\end{equation}
is non-negative. Here, $\theta_{1,i}$ reflects the user's preference for faster charging, and $\theta_{2,i}$ captures their sensitivity to price markup.
To introduce meaningful behavioral heterogeneity, we draw both parameters independently from uniform distributions, $\theta_{1,i} \sim \mathcal{U}(0.5,2.0)$ and $\theta_{2,i} \sim \mathcal{U}(0.5,2.0)$.
This setting allows for a wide range of user types. Lower values of $\theta_{1,i}$ correspond to users who are willing to accept longer charging durations, whereas higher values indicate a stronger preference for faster charging. The parameter $\theta_{2,i}$ captures heterogeneity in users’ sensitivity to price markups. Together, these parameter settings introduce realistic heterogeneity in user acceptance behavior while maintaining comparability between charging-rate preference and price sensitivity.

\begin{figure}
  \centering
 \includegraphics[width=.8\textwidth]{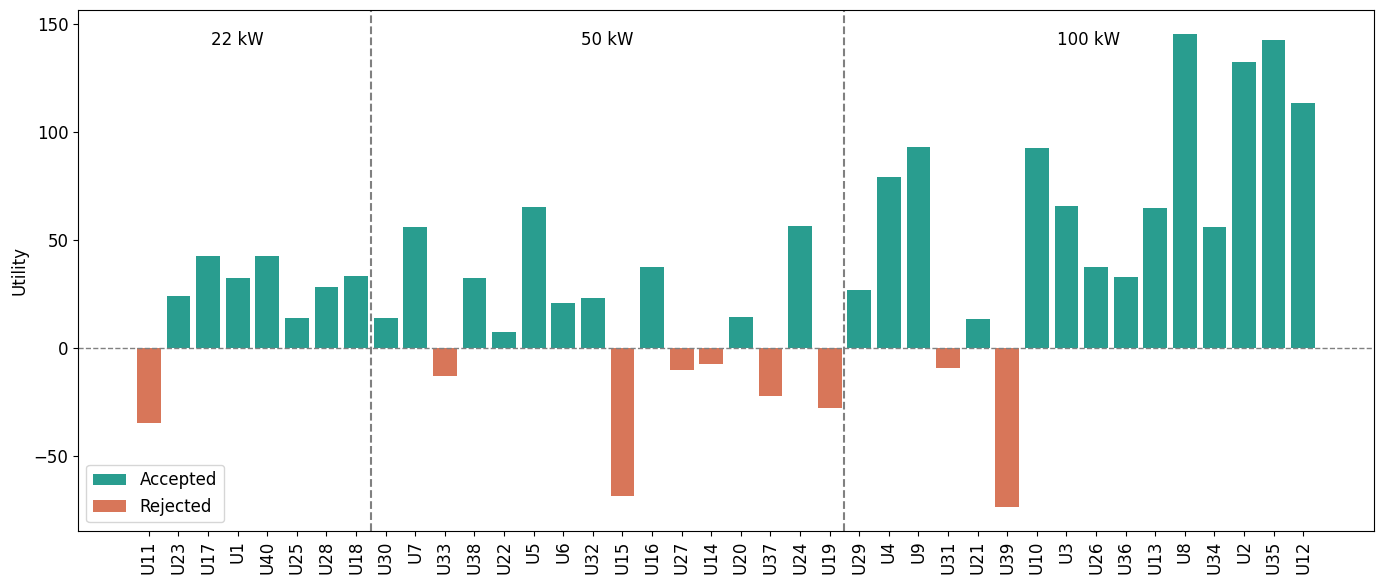}
  \caption{User Acceptance Based on Utility Function}\label{fig:user_acceptance}
\end{figure}

Among all allocations accepted by the operator, we further evaluate whether individual users would ultimately confirm or reject the assigned plans based on their personalized utility scores.
We next relate the user-level utility outcomes in Figure \ref{fig:user_acceptance} to the bid and counter-price patterns shown in Figure \ref{fig:40users_48slots_price_wide}.
Figure \ref{fig:40users_48slots_price_wide} shows that for the 22~kW option, only one user experiences a price markup, while all remaining users are accepted at their bid prices. Correspondingly, Figure \ref{fig:user_acceptance} indicates that users without price markups achieve positive utilities, whereas the user facing a markup exhibits a negative utility and rejects the allocation.
At 50 kW, all users are subject to price markups. Despite this, Figure \ref{fig:user_acceptance} shows that the majority of users still obtain positive utilities and accept the assigned plans, confirming that markups can be acceptable when accompanied by higher charging rates.
For the 100 kW option, most users face price markups, while both positive and a small number of negative utilities are observed. Notably, the positive utility values are generally higher than those observed at lower charging rates, indicating that higher charging rates can generate larger net utility gains for many users despite the presence of markups.
Overall, these results confirm that price markups induced by the operator’s counter-offers can lead to negative user utility and rejection, while also showing that higher charging rates can offset markups for a substantial fraction of users.

\subsection{Discussion and operational insights}

{The results provide several operational insights for the design and management of public EV charging systems, particularly regarding how pricing flexibility, user participation, and infrastructure configuration jointly shape overall system performance.}

{First, the price preservation parameter $\gamma$ highlights a fundamental trade-off between pricing flexibility and perceived fairness. Higher values of $\gamma$ enforce stronger price integrity and improve user satisfaction by restricting the operator’s ability to increase prices above user bids. However, this constraint also limits the operator’s flexibility to adjust prices and allocate high-power charging resources profitably. This tension is particularly pronounced for high-power chargers (e.g., 100 kW), whose economic viability depends strongly on sufficient pricing flexibility. From an operational perspective, these findings suggest that allowing moderate pricing flexibility can help operators better coordinate charging demand and improve infrastructure utilization. From a policy perspective, overly strict price preservation rules may limit the effectiveness of demand coordination mechanisms in public charging systems.}

{From a behavioral perspective, user acceptance is governed more by absolute pricing thresholds than by relative bid rankings. Users submitting very low bids are frequently rejected or subject to markups, while moderate bids often remain acceptable within the operator’s pricing margins. Interestingly, higher bids do not always guarantee access to premium charging if they fail to meet the operator’s internal thresholds. This suggests that transparent pricing feedback and learning mechanisms may play an important role in repeated platform interactions. Providing users with clearer information about acceptance thresholds or expected prices could improve bidding strategies and reduce repeated rejections, thereby improving overall charging system efficiency.}

{The distributional outcomes further illustrate how pricing flexibility shapes access to charging resources. Gini coefficients increase with higher $\gamma$ values and higher charging rates, indicating that greater pricing flexibility can lead to wider dispersion among accepted users. However, inequality remains relatively low across scenarios, with Gini coefficients consistently below 0.1. This suggests that the proposed framework can balance efficiency improvements with a relatively equitable allocation of charging opportunities. This result suggests that coordinated pricing mechanisms can improve infrastructure utilization without substantially increasing inequality among users.}

{Finally, heterogeneous charger configurations can substantially improve system inclusiveness and operational efficiency. Combining 22 kW and 50 kW chargers enables better alignment between user demand and charging capabilities, particularly under stricter price preservation settings. In contrast, single-rate systems tend to exclude low-bidding or time-constrained users under tight constraints. These findings suggest that infrastructure design and pricing rules should be jointly considered when planning public charging systems, rather than treated as independent decisions. More broadly, the results indicate that pricing mechanisms, charging scheduling, and infrastructure configuration should be jointly considered when designing public charging systems to achieve both efficiency and fairness.}

\section{Conclusion and future work} \label{Section5}

{This paper proposes a user-driven pricing and scheduling framework for public EV charging that integrates user bidding with operator-side optimization under capacity constraints. By allowing users to express their charging preferences through bids while enabling operators to coordinate prices and charging slots centrally, the framework captures the interaction between user incentives, pricing flexibility, and limited charging capacity. The proposed approach provides a tractable mechanism for coordinating heterogeneous charging demand while balancing pricing fairness and system efficiency.}

{The numerical results highlight several key insights regarding the design of coordinated EV charging systems. First, pricing flexibility plays a central role in balancing perceived fairness and operational efficiency. While stronger price preservation improves user acceptance, it also reduces the operator’s ability to allocate charging resources efficiently, particularly for high-power chargers. Second, user acceptance is shaped more by absolute pricing thresholds than by relative bid rankings, suggesting that transparent pricing feedback can improve user participation and bidding behavior. Third, charger heterogeneity naturally segments users and improves the alignment between charging demand and infrastructure capabilities. Together, these findings suggest that pricing mechanisms, charging scheduling, and infrastructure configuration should be jointly considered when designing public charging systems.}

{The current framework assumes that bids are collected before scheduling, whereas real-world charging systems operate with asynchronous user arrivals and require rapid responses. Future research could therefore investigate rolling or event-driven coordination mechanisms that update prices and allocations as new bids arrive.
On the behavioral side, incorporating richer user response models, such as probabilistic acceptance or strategic bidding behavior, could further improve realism. On the system side, integrating grid constraints, real-time electricity prices, and demand uncertainty would enhance the robustness of charging coordination.
Finally, extending the framework to multi-station environments would allow the study of network-level charging systems where users choose among multiple stations. Such extensions would introduce additional challenges, including station choice behavior, cross-station allocation mechanisms, and potential competition between operators.}

\section{Acknowledgements}

The authors gratefully acknowledge the support of Chalmers University of Technology, Area of Advanced Transport, through the project User-Driven Market Mechanisms for Public EV Charging. The authors also acknowledge the support of the ERGODIC project (F-DUT-2022-0078), funded by the European Commission and Vinnova; the E-Laas project (F-ENUAC-2022-0003), supported by the European Commission and the Swedish Energy Agency; and the Game of Own project (2024-01053), funded by the Swedish Energy Agency.

\appendix
\section*{Appendix A} \label{AppendixA}

\setcounter{table}{0}  
\renewcommand{\thetable}{A\arabic{table}}

Table \ref{parameters} summarizes the sets, decision variables, and parameters used in the proposed optimization model.

\begin{table}[width=.9\linewidth,cols=4,pos=h]
\centering
\caption{Sets and parameters}\label{parameters}
		\begin{tabular*}{\tblwidth}{@{} LL@{} }
\toprule
\textbf{Sets:}&\\
$ \mathcal{I} $ & Set of all users, $i \in \mathcal{I}$ \\
$ \mathcal{T}_i $ & Set of acceptable time slots for user $ i $, $t \in \mathcal{T}_i$ \\
$ \mathcal{R} $ & Set of charging rates, $r \in \mathcal{R}$ \\
\midrule
\textbf{Decision Variables:} &\\
$a_i$ & 1 if user $i$ is accepted to be served, 0 otherwise.\\
$x^{\text{B}}_{i,t,r}$ & 1 if the user’s original bid is accepted, 0 otherwise. \\
$x^{\text{N}}_{i,t,r}$ & 1 if a new price is offered by the operator, 0 otherwise. \\
$p^{\text{N}}_{i,t,r}$ & New price proposed by the operator for $(i,t,r)$. \\
$p^{\text{F}}_{i,t,r}$ & Final price: $x^{\text{B}}_{i,t,r} p^{\text{B}}_i + x^{\text{N}}_{i,t,r} p^{\text{N}}_{i,t,r}$. \\
$y_{i,r}$ & 1 if user $i$ is assigned to charging rate $r$ by the operator, 0 otherwise. \\
$z_{i,t,r}$ & 1 if user $i$ is assigned to time slot $t$ at rate $r$, 0 otherwise. \\
$z^{\text{s}}_{i,t}, z^{\text{e}}_{i,t}$ & 1 if slot $t$ is the start/end of user $i$’s charging window. \\
$q_{i,t}$ & Actual energy delivered to user $i$ in slot $t$. \\
$p^{\text{R}}_{t,r}$ & Reference price for each slot and rate. \\
\midrule
\textbf{Parameters:} &\\
$p^{\text{B}}_i$ & Bid price per kWh from user $i$. \\
$A_{i,t}$ & Equals 1 if time slot $t$ is accepted by user $i$ ($t \in \mathcal{T}_i$), and 0 otherwise. \\
$Q_i^{\min}, Q_i^{\max}$ & User $i$’s min and max charging demand. \\
$E_r$ & Energy delivered per time slot for rate $r$: $E_r = R_r \cdot \tau / 60$. \\
$M_r$ & Number of chargers available at rate $r$. \\
$R_r$ & Charging rate (kW) for level $r \in \{1,2,3\}$. \\
$\tau$ & Length of each time slot (slot duration).\\
$p^{\text{H}}$ & Highest charging price, serving as the big M.\\
$c_{r,t}$ & Unit cost of charging at rate $r$ during time slot $t$. \\
$C_t$ & System-wide slot capacity in kWh at time $t$. \\
$\delta_r$ & Maximum allowed markup over user bid at rate $r$. \\
$\alpha$ & Maximum price-to-cost ratio. \\
$\varepsilon$ & Minimum profit margin. \\
$\gamma$ & 	Slot-level bid price preservation ratio. \\
\bottomrule
		\end{tabular*}
\end{table}


\bibliographystyle{cas-model2-names}

\bibliography{reference}


\end{document}